\documentclass[a4paper,11pt]{article}

\usepackage{amsfonts,amssymb,amsmath,amsthm,latexsym,bbm}
\usepackage{epsf,epsfig}
\usepackage{xcolor,colortbl,color}
\usepackage{graphicx,graphics}
\usepackage{caption}

\usepackage[english]{babel}
\makeatletter \@addtoreset{equation}{section} \makeatother
\makeatletter \@addtoreset{enunciato}{section} \makeatother
\newcounter{enunciato}[section]

\newtheorem{ittheorem}{Theorem}
\newtheorem{itlemma}{Lemma}
\newtheorem{itproposition}{Proposition}
\newtheorem{itdefinition}{Definition}
\newtheorem{itremark}{Remark}
\newtheorem{itclaim}{Claim}
\newtheorem{itfact}{Fact}
\newtheorem{itconjecture}{Conjecture}
\newtheorem{itcorollary}{Corollary}

\newenvironment{theorem}{\addtocounter{enunciato}{1}
\begin{ittheorem}}{\end{ittheorem}}
\newenvironment{lemma}{\addtocounter{enunciato}{1}
\begin{itlemma}}{\end{itlemma}}
\newenvironment{proposition}{\addtocounter{enunciato}{1}
\begin{itproposition}}{\end{itproposition}}
\newenvironment{definition}{\addtocounter{enunciato}{1}
\begin{itdefinition}}{\end{itdefinition}}
\newenvironment{remark}{\addtocounter{enunciato}{1}
\begin{itremark}}{\end{itremark}}

\newenvironment{conjecture}{\addtocounter{enunciato}{1}
\begin{itconjecture}}{\end{itconjecture}}
\newenvironment{corollary}{\addtocounter{enunciato}{1}
\begin{itcorollary}}{\end{itcorollary}}

\newcommand{\be}[1]{\begin{equation}\label{#1}}
\newcommand{\ee}{\end{equation}}
\newcommand{\bl}[1]{\begin{lemma}\label{#1}}
\newcommand{\el}{\end{lemma}}
\newcommand{\br}[1]{\begin{remark}\label{#1}}
\newcommand{\er}{\end{remark}}
\newcommand{\bt}[1]{\begin{theorem}\label{#1}}
\newcommand{\et}{\end{theorem}}
\newcommand{\bd}[1]{\begin{definition}\label{#1}}
\newcommand{\ed}{\end{definition}}
\newcommand{\bp}[1]{\begin{proposition}\label{#1}}
\newcommand{\ep}{\end{proposition}}
\newcommand{\bc}[1]{\begin{corollary}\label{#1}}

\newcommand{\bcj}[1]{\begin{conjecture}\label{#1}}
\newcommand{\ecj}{\end{conjecture}}

\newcommand{\R}{\ensuremath{\mathbb{R}}}
\newcommand{\var}{\ensuremath{\mathbb{V}\mathrm{ar}}}
\newcommand{\Z}{\ensuremath{\mathbb{Z}}}
\newcommand{\N}{\ensuremath{\mathbb{N}}}
\newcommand{\PP}{\ensuremath{\mathbb{P}}}
\newcommand{\QQ}{\ensuremath{\mathbb{Q}}}
\newcommand{\E}{\ensuremath{\mathbb{E}}}
\newcommand{\abs}[1]{\ensuremath{\left|#1\right|}}
\newcommand{\dd}{\ensuremath{\mathrm{d}}}

\newcommand{\argmax}{\ensuremath{\text{argmax}}}

\newcommand{\cM}{\mathcal{M}} 
\newcommand{\cN}{\mathcal{N}}
\newcommand{\cE}{\mathcal{E}}

\begin{document}
\title{Random walks in cooling random environments}

\author{\renewcommand{\thefootnote}{\arabic{footnote}}
L.\ Avena\,\footnotemark[1]\,\,, F.\ den Hollander\,\footnotemark[1]}

\date{\today}

\footnotetext[1]{Mathematical Institute, Leiden University, 
P.O.\ Box 9512, 2300 RA Leiden, The Netherlands}

\maketitle

\begin{abstract}
We propose a model of a one-dimensional random walk in dynamic random 
environment that interpolates between two classical settings: (I) the random 
environment is sampled at time zero only; (II) the random environment is 
resampled at every unit of time. In our model the random environment is resampled 
along an increasing sequence of deterministic times. We consider the annealed 
version of the model, and look at three growth regimes for the resampling times: 
(R1) linear; (R2) polynomial; (R3) exponential. We prove weak laws of large 
numbers and central limit theorems. We list some open problems and conjecture 
the presence of a crossover for the scaling behaviour in regimes (R2) and (R3).

\\
{\it Acknowledgment:} The research in this paper was supported through ERC Advanced 
Grant VARIS--267356 and NWO Gravitation Grant NETWORKS-024.002.003. The authors 
are grateful to David Stahl for his input at an early stage of the project, and to Zhan Shi for 
help with the argument in Appendix~\ref{appLp}. Thanks also to Yuki Chino and Conrado 
da Costa for comments on a draft of the paper.
\end{abstract}

\newpage

%%%%%%%%% SECTION 1 %%%%%%%%%%%%%%

\section{Introduction, model, main theorems and discussion}
\label{s1}

%%%%%%%%%

\subsection{Background and outline}
\label{ssintro}

Models for particles moving in media with impurities pose many challenges. In 
mathematical terms, the medium with impurities is represented by a random 
environment on a lattice and the particle motion is represented by a random 
walk on this lattice with transition probabilities that are determined by the environment.

\medskip\noindent
{\bf Static Random Environment.}
A model that has been studied extensively in the literature is that of \emph{Random 
Walk in Random Environment} (RWRE), where the random environment is \emph{static} 
(see Zeitouni~\cite{Z04} for an overview). In one dimension this model exhibits striking 
features that make the presence of a random environment particularly interesting. 
Namely, there are regions in the lattice where the random walk remains \emph{trapped} 
for a very long time. The presence of these traps leads to a local slow down of the random 
walk in comparison to a homogeneous random walk. This, in turn, is responsible for a 
non-trivial limiting speed, as well as for anomalous scaling behaviour (see Section~\ref{ssRWRE}
below). For instance, under proper assumptions on the random environment, the 
random walk can be transient yet sub-ballistic in time, or it can be non-diffusive in time 
with non-Gaussian fluctuations. To derive such results and to characterise associated 
limit distributions, a key approach has been to represent the environment by a potential 
function: a deep valley in the potential function corresponds to a region in the lattice 
where the random walk gets trapped for a very long time. 

\medskip\noindent
{\bf Dynamic Random Environment.} 
If, instead, we consider a random walk in a random environment that itself evolves over 
time, according to a prescribed \emph{dynamic} rule, then the random environment is 
still inhomogeneous, but the dynamics dissolves existing traps and creates new traps. 
Depending on the choice of the dynamics, the random walk behaviour can be similar 
to that in the static model, i.e., show some form of \emph{localisation}, or it can be similar 
to that of a homogeneous random walk, in which case we speak of \emph{homogenisation}. 
The simplest model of a dynamic random environment is given by an i.i.d.\ field of spatial
random variables that is resampled in an i.i.d.\ fashion after every step of the random walk. 
This model has been studied in several papers. Clearly, under the so-called annealed measure
homogenisation occurs: the independence of the random environment in space and time
causes the random walk to behave like a homogeneous random walk, for which a standard 
law of large numbers and a standard central limit theorem hold.

\medskip\noindent
{\bf Cooling Random Environment.}
In the present paper we introduce a new model, which we call \emph{Random Walk in 
Cooling Random Environment} (RWCRE). This model interpolates between the two 
settings mentioned above: start at time zero with an i.i.d.\ random environment and 
resample it along an increasing sequence of deterministic times (the name ``cooling" 
is chosen here because the static model is sometimes called ``frozen"). If the resampling 
times increase rapidly enough, then we expect to see a behaviour close to that of the 
static model, via some form of \emph{localisation}. Conversely, if the resampling 
times increase slowly enough, then we expect \emph{homogenisation} to be dominant. 
The main goal of our paper is to start making this rough heuristics precise by proving 
a few basic results for RWCRE in a few cooling regimes. From a mathematical point of 
view, the analysis of RWCRE reduces to the study of sums of independent random 
variables distributed according to the static model, viewed at successive time scales. 
In Section~\ref{discussion} we will point out what type of results for the static model, 
\emph{currently still unavailable}, would be needed to pursue a more detailed analysis of 
RWCRE.   

\noindent
{\bf Outline.}
In Section~\ref{ssRWRE} we recall some basic facts about the one-dimensional RWRE model. 
In Section~\ref{ssRWCRE} we define our model with a \emph{cooling} random environment 
and introduce the three cooling regimes we are considering. In Section~\ref{results} we state 
three theorems and make a remark about more general cooling regimes. In Section~\ref{discussion} 
we mention a few open problems and state a conjecture. Sections~\ref{s2}--\ref{s4} are devoted 
to the proofs of our theorems. Appendix~\ref{app:Toep} proves a variant of a Toeplitz lemma that
is needed along the way. Appendix~\ref{app:LL} recalls the central limit theorem for sums of 
independent random variables. Appendix~\ref{appLp} strengthens a well-known convergence 
in distribution property of recurrent RWRE, which is needed for one of our theorems.

%%%%%%%%%

\subsection{RWRE}
\label{ssRWRE}

Throughout the paper we use the notation $\N_0 = \N \cup \{0\}$ with $\N = \{1,2,\dots\}$. 
The classical one-dimensional random walk in random environment (RWRE) is defined 
as follows. Let $\omega=\{\omega(x)\colon\,x\in\Z\}$ be an i.i.d.\ sequence with probability 
distribution 
\begin{equation}
\mu = \alpha^{\Z}
\end{equation}
for some probability distribution $\alpha$ on $(0,1)$. The random walk in the \emph{space} 
environment $\omega$ is the Markov process $Z=(Z_n)_{n\in\N_0}$ starting at $Z_0=0$ with transition 
probabilities
\be{RWRE}
P^\omega(Z_{n+1}=x+e \mid Z_n=x)
= \left\{\begin{array}{ll}
\omega(x), &\mbox{ if } e=1,\\
1-\omega(x), &\mbox{ if } e=-1,
\end{array}
\right. \qquad n\in\N_0.
\ee
The properties of $Z$ are well understood, both under the \emph{quenched} law
$P^\omega(\cdot)$ and the \emph{annealed} law 
\be{RWREann}
\PP_\mu(\cdot) = \int_{(0,1)^\Z} P^\omega(\cdot)\,\mu(\dd\omega). 
\ee
Let $\langle\cdot\rangle$ denote expectation w.r.t.\ $\alpha$. Abbreviate
\be{rhodef}
\rho(0) = \frac{1-\omega(0)}{\omega(0)}.
\ee
Without loss of generality we make the following assumption on $\alpha$:
\begin{equation}
\label{Nested}
\langle \log\rho(0) \rangle \leq 0.
\end{equation}
This assumption guarantees that $Z$ has a preference to move to the right. In what
follows we recall some key results for RWRE on $\Z$ . For a general overview, we 
refer the reader to Zeitouni~\cite{Z04}. 

The following result due to Solomon \cite{S75} characterises recurrence versus transience 
and asymptotic speed.

\bp{Solomon} {\bf [Recurrence, transience and speed RWRE]}
Let $\alpha$ be any distribution on $(0,1)$ satisfying~\eqref{Nested}. Then the 
following hold:  
\begin{itemize} 
\item
$Z$ is recurrent if and only if $\langle\log\rho(0)\rangle = 0$.
\item
If $\langle\log\rho(0)\rangle < 0$ and $\langle\rho(0)\rangle \geq 1$, then $Z$ is transient to the 
right with zero speed:
\begin{equation}
\lim_{n\to\infty} \frac{Z_n}{n} = 0 \quad \PP_\mu\text{-a.s.}
\end{equation}
\item
If $\langle\log\rho(0)\rangle < 0$ and $\langle\rho(0)\rangle <1$, then $Z$ is transient to the right 
with positive speed: 
\begin{equation}
\label{vformula}
\lim_{n\to\infty} \frac{Z_n }{n} = v_\mu = \frac{1-\langle\rho(0)\rangle}{1+\langle\rho(0)\rangle} > 0 
\quad \PP_\mu\text{-a.s.}
\end{equation}
\end{itemize}
\ep

The scaling limits in the different regimes have been studied in a number of papers,
both under the quenched and the annealed law. While the results are the same for 
the law of large numbers, they are in general different for the scaling limits. Under 
the quenched law only partial results are available (see Peterson~\cite{Peterson}
for a summary of what is known). For this reason we are forced to restrict ourselves 
to the annealed law.

In the recurrent case the proper scaling was identified by Sinai~\cite{S82} and the limit law 
by Kesten~\cite{K86}. The next proposition summarises their results.

\bp{RecScaling} {\bf [Scaling limit RWRE: recurrent case]}
Let $\alpha$ be any probability distribution on $(0,1)$ satisfying $\langle\log\rho(0)\rangle 
= 0$ and $\sigma_\mu^2 = \langle\log^2\rho(0)\rangle \in (0,\infty)$. Then, under the annealed 
law $\PP_\mu$, the sequence of random variables 
\be{Sinai}
\frac{Z_n}{\sigma_\mu^2\log^2 n}, \quad n\in\N,
\ee 
converges in distribution to a random variable $V$ on $\R$ that is independent of $\alpha$.
The law of $V$ has a density $p(x)$, $x\in\R$, with respect to the Lebesgue measure that 
is given by 
\be{densityofV}
p(x) = \frac{2}{\pi} \sum_{k\in\N_0} \frac{(-1)^k}{2k+1}
\exp\left[-\frac{(2k+1)^2\pi^2}{8}|x| \right], \qquad x \in \R.
\ee
\ep

\noindent
For later use we need to show that the sequence in~\eqref{Sinai} converges to $V$ in $L^p$ 
for every $p>0$. This is done in Appendix~\ref{appLp}. Note that the law of $V$ is symmetric 
with variance $\sigma^2_V \in (0,\infty)$.

The scaling in the (annealed) transient case was first studied by Kesten, Kozlov and 
Spitzer~\cite{KKS75}. In order to state their results, we need some more notation. 
Given $s,b>0$, denote by $L_{s,b}$ the $s$-stable distribution with scaling parameter 
$b$, centred at $0$ and totally skewed to the right. In formulas, $L_{s,b}$ lives on $\R$ 
and is identified by its characteristic function
\be{stableDensity}
\hat{L}_{s,b}(u) = \int_\R e^{i ux } L_{s,b}(dx)
=\exp\left[-b |u|^s \left(1-i\frac{u}{|u|}g_s(u)\right)\right], \qquad u \in \R,
\ee
with 
\begin{equation}
g_s(u) = \left\{\begin{array}{ll}
\tan(\frac{s\pi}{2}),   &s \neq 1,\\
\frac{2}{\pi}\log{|u|}, &s =1.
\end{array}
\right. \qquad u \in \R.
\end{equation}
Write $\phi(x) = \frac{1}{\sqrt{2\pi}} \int_{-\infty}^x e^{-y^2/2} dy$, $x \in \R$, to denote 
the standard normal distribution.

\bp{TransScaling} {\bf [Scaling limit RWRE: transient case]}
Let $\alpha$ be any probability distribution on $(0,1)$ satisfying $\langle\log\rho(0)\rangle 
< 0$ such that the support of the distribution of $\log\rho(0)$ is non-lattice. Let $s \in (0,\infty)$ 
be the unique root of the equation
\begin{equation}
\label{sparamdef}
\langle\rho(0)^s\rangle = 1,
\end{equation} 
and suppose that $\langle\rho(0)^s\log\rho(0)\rangle<\infty$. Then, under the annealed 
law $\PP_\mu$, the following hold:
\begin{itemize}
\item 
If $s\in(0,1)$, then there exists a $b>0$ such that 
\be{inverseStable}
\lim_{n\to\infty} \PP_{\mu}\left( \frac{Z_n}{n^{s}} \leq x\right)
= [1-L_{s,b}(x^{-1/s})]\,{\mathbbm 1}_{\{x>0\}}.
\ee 
\item 
If $s=1$, then there exist $b>0$ and $\{\delta_\alpha(n)\}_{n\in\N}$, satisfying 
$\delta_\alpha(n)=[1+o(1)]\,n/b\log n$ as $n\to\infty$, such that 
\be{Cauchy}
\lim_{n\to\infty} \PP_{\mu}\left(\frac{Z_n-\delta_\alpha(n)}{n/\log^2 n} \leq x\right)
= 1-L_{1,b}(-b^2 x), \qquad x \in \R.
\ee 
\item 
If $s\in(1,2)$, then there exist $b>0$ and $c=c(b)>0$ such that 
\be{StableLeftSkew}
\lim_{n\to\infty} \PP_{\mu}\left(\frac{Z_n-v_{\mu}n}{b n^{1/s}}\leq x\right)
= 1-L_{s,c}(-x), \qquad x \in \R.
\ee 
\item
If $s=2$, then there exists a $b>0$ such that 
\be{AlmostCLT}
\lim_{n\to\infty} \PP_{\mu}\left(\frac{Z_n-v_{\mu}n}{ b\sqrt{n\log n}} \leq x\right)
= \phi(x), \qquad x \in \R.
\ee 
\item
If $s\in(2,\infty)$, then there exists a $b>0$ such that 
\be{CLT}
\lim_{n\to\infty} \PP_{\mu}\left(\frac{Z_n-v_{\mu}n}{ b\sqrt{n}} \leq x\right)
= \phi(x), \qquad x \in \R.
\ee 
\end{itemize}  
\ep

\noindent
In~\eqref{Cauchy} and~\eqref{StableLeftSkew}, the limiting laws are stable 
laws that are totally skewed to the left, i.e., their characteristic function is as in 
\eqref{stableDensity} but with a $+$ sign in the term with the imaginary factor 
$i$ in the exponential. In~\eqref{inverseStable}, the limiting law is an inverse 
stable law, sometimes referred to as the Mittag-Leffler distribution. 

%%%%%%%%%%%%%%

\subsection{RWCRE}
\label{ssRWCRE}

\begin{figure}[htbp]
\vspace{-.5cm}
\begin{center}
\setlength{\unitlength}{1cm}
\begin{picture}(13,2)(-.5,-1)
\put(.1,.1){\line(1,0){3}}
\put(4.1,.1){\line(1,0){4}}
\put(9.1,.1){\line(1,0){4}}
{\thicklines
\qbezier(13.1,-.1)(13.1,0)(13.1,.25)
}
\qbezier[10](3.4,.1)(3.6,.1)(3.8,.1)
\qbezier[10](8.4,.1)(8.6,.1)(8.8,.1)
\put(0,0){$\bullet$}
\put(2,0){$\bullet$}
\put(5,0){$\bullet$}
\put(7,0){$\bullet$}
\put(10,0){$\bullet$}
\put(-.3,-.5){$\tau(0)$}
\put(1.7,-.5){$\tau(1)$}
\put(4.3,-.5){$\tau(k-1)$}
\put(6.7,-.5){$\tau(k)$}
\put(9.5,-.5){$\tau(k(n))$}
\put(13,-.5){$n$}
\end{picture}
\end{center}
\vspace{-.5cm}
\caption{Resampling times $\tau(k)$, $0\leq k \leq k(n)$, prior to time $n$.}
\label{fig:resamp1}
\vspace{.3cm}
\end{figure}
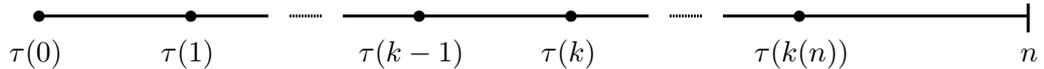

In the present paper we look at a model where $\omega$ is updated along a growing sequence 
of deterministic times. To that end, let $\tau\colon\,\N_0 \to \N_0$ be a strictly increasing map 
such that $\tau(0)=0$ and $\tau(k) \geq k$ for $k\in\N$ (see Fig.~\ref{fig:resamp1}). Define a 
sequence of random environments $\Omega = (\omega_n)_{n\in\N_0}$ as follows: 
\begin{itemize}
\item[$\blacktriangleright$]
At each time $\tau(k)$, $k\in\N_0$, the environment $\omega_{\tau(k)}$ is freshly resampled from 
$\mu=\alpha^\Z$ and does not change during the time interval $[\tau(k),\tau(k+1))$.
\end{itemize} 

%%%%%%%%%%%%%%%%%%%%%%%%%%%%%%%%%%%%%%%%%
%%%%%%%%%%%%%%%%%%%%%%%%%%%%%%%%%%%%%%%%%

The random walk in the \emph{space-time} environment $\Omega$ is the Markov process 
$X=(X_n)_{n\in\N_0}$ starting at $X_0=0$ with transition probabilities
\be{SignStSpeed}
P^\Omega(X_{n+1}=x+e \mid X_n=x)
= \left\{\begin{array}{ll}
\omega_n(x), &\mbox{ if } e=1,\\
1-\omega_n(x), &\mbox{ if } e=-1,
\end{array}
\right. \qquad n \in \N_0.
\ee
We call $X$ the \emph{random walk in cooling random environment} (RWCRE) with 
\emph{resampling rule} $\alpha$ and \emph{cooling rule} $\tau$. Our goal will be to 
investigate the behavior of $X$ under the \emph{annealed law}
\begin{equation}
\label{RWCREann}
\PP_\QQ(\cdot) = \int_{(\R^\Z)^{\N_0}} P^\Omega(\cdot)\,\QQ(\dd\Omega),
\end{equation} 
where $\QQ=\QQ_{\alpha,\tau}$ denotes the law of $\Omega$.

Note that if $\tau(1) = \infty$, then RWCRE reduces to RWRE, i.e., $X$ has the 
same distribution as $Z$. On the other hand, if $\tau(k)=k$, $k\in\N$, then $\omega_n$ 
is freshly sampled from $\mu$ for all $n\in\N$ and RWCRE reduces to what is often referred 
to as random walk in an i.i.d.\ space-time random environment. Under the annealed law, 
the latter is a homogeneous random walk and is trivial. Under the quenched law it is 
non-trivial and has been investigated in a series of papers, e.g.\ Boldrighini, Minlos, 
Pellegrinotti and Zhizhina~\cite{BMPseries}, and Rassoul-Agha and Seppalainen~\cite{RAS05}. 
For any other choice of $\tau$ the model has not been considered before. 

In this paper we will focus on \emph{three different growth regimes} for $\tau(k)$ as $k\to\infty$:
\begin{enumerate}
\item[(R1)]
\emph{No cooling:}
$\tau(k) \sim Ak$ for some $A \in (1,\infty)$.
\item[(R2)] 
\emph{Slow cooling:}
$\tau(k) \sim Bk^{\beta}$ for some $B \in (0,\infty)$ and $\beta \in (1,\infty)$.
\item[(R3)] 
\emph{Fast cooling:} 
$\log \tau(k) \sim Ck$ for some $C \in (0,\infty)$.
\end{enumerate}

Let 
\be{TimeIncr}
T_k = \tau(k)-\tau(k-1), \quad k\in\N,
\ee
be the increments of the resampling times of the random environment. We assume
that $k \mapsto T_k$ is sufficiently regular so that 
\be{Tnreg}
\begin{array}{ll}
T_k \sim \beta B k^{\beta-1} &\text{ in regime (R2)},\\
\log T_k \sim Ck &\text{ in regime (R3)}. 
\end{array}
\ee
For instance, in regime (R2) this regularity holds as soon as $k \mapsto T_k$ is 
ultimately non-decreasing (Bingham, Goldie and Teugels~\cite[Sections 1.2 and 
1.4]{BGT87}).

%%%%%%%%%%%%%%%%%%%%%%%%%%%%%%%%%%%%%%%%%
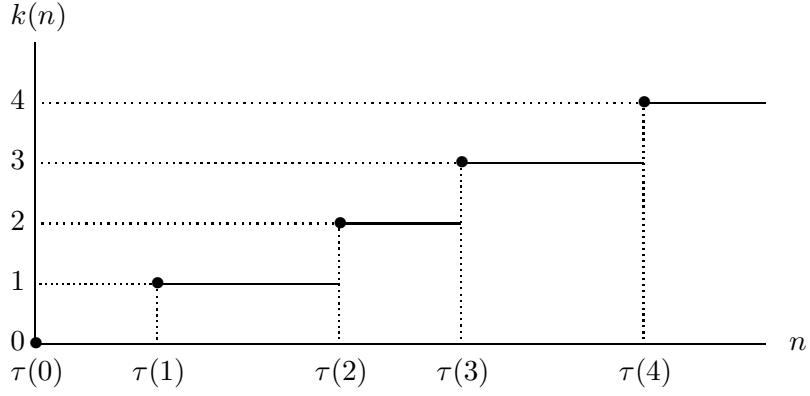
\begin{figure}[htbp]
\vspace{1.2cm}
\begin{center}
\setlength{\unitlength}{0.8cm}
\begin{picture}(13,5)(-.5,-.2)
\put(.1,.1){\line(1,0){12}}
\put(.1,.1){\line(0,1){5}}
%%%
\put(.1,.1){\line(1,0){2}}
\put(2.1,1.1){\line(1,0){3}}
\put(5.1,2.1){\line(1,0){2}}
\put(7.1,3.1){\line(1,0){3}}
\put(10.1,4.1){\line(1,0){2}}
\qbezier[20](.1,1.1)(1,1.1)(2,1.1)
\qbezier[40](.1,2.1)(2.5,2.1)(5,2.1)
\qbezier[60](.1,3.1)(3.5,3.1)(7,3.1)
\qbezier[80](.1,4.1)(5,4.1)(10,4.1)

\qbezier[10](2.1,.1)(2.1,.55)(2.1,1.1)
\qbezier[20](5.1,.1)(5.1,1.05)(5.1,2.1)
\qbezier[30](7.1,.1)(7.1,1.65)(7.1,3.1)
\qbezier[40](10.1,.1)(10.1,2.05)(10.1,4.1)

%%%
\put(12.5,0){$n$}
\put(-.3,5.4){$k(n)$} 
\put(-.3,-.5){$\tau(0)$}
\put(1.7,-.5){$\tau(1)$}
\put(4.7,-.5){$\tau(2)$}
\put(6.7,-.5){$\tau(3)$}
\put(9.7,-.5){$\tau(4)$}
\put(-.3,0){$0$}
\put(-.3,1){$1$}
\put(-.3,2){$2$}
\put(-.3,3){$3$}
\put(-.3,4){$4$}
\put(0,0){$\bullet$}
\put(2,1){$\bullet$}
\put(5,2){$\bullet$}
\put(7,3){$\bullet$}
\put(10,4){$\bullet$}
\end{picture}
\end{center}
\vspace{.2cm}
\caption{Plot of $n \mapsto k(n)$. }
\label{fig:resampkn}
\end{figure}
%%%%%%%%%%%%%%%%%%%%%%%%%%%%%%%%%%%%%%%%%

Let
\be{SpaceIncr}
Y_k = X_{\tau(k)}-X_{\tau(k-1)}, \quad k\in\N,
\ee
be the increments of the random walk between the resampling times. Our 
starting point will be the relation
\be{XRepresentation}
X_n = \sum_{k=1}^{k(n)} Y_k + \bar Y^n,
\ee
where (see Fig.~\ref{fig:resampkn})
\be{taudefs}
k(n) = \max \{k \in \N\colon\,\tau(k) \leq n\}
\ee
is the last resampling prior to time $n$, $Y_k$ is distributed as $Z_{T_k}$ in 
environment $\omega_{\tau(k-1)}$, while $\bar  Y^n$ is a boundary term that 
is distributed as $Z_{\bar  T^n}$ in environment $\omega_{\tau(k(n))}$ with
\be{barTdef}
\bar T^n = n - \tau(k(n))
\ee
the remainder time after the last resampling. Note that all terms in \eqref{XRepresentation} 
are independent.

%%%%%%%%%%%

\subsection{Main theorems}
\label{results}

We are now ready to state our results under the annealed measure $\PP_\QQ$ 
in the three cooling regimes (R1)--(R3).

\medskip\noindent
$\bullet$ {\bf No cooling.}
Regime (R1) in essence corresponds to the situation where the increments of the resampling 
times do not diverge, i.e., $\lim_{k\to\infty} T_k \neq \infty$. To analyze this regime, we assume 
that the empirical measure of the increments of the resampling times
\be{empiricaldistribution}
L_n = \frac{1}{k(n)} \sum_{k=1}^{k(n)} \delta_{\tau(k)-\tau(k-1)}
\ee
has a non-degenerate $L^1$-limit, i.e., 
\be{weakconvergence}
\lim_{n\to\infty} \sum_{\ell\in\N} \ell\,|L_n(\ell)-\nu(\ell)| = 0
\ee
for some $\nu \in \cM_1(\N)$. Since $\sum_{\ell\in\N} \ell L_n(\ell) = \tau(k(n))/k(n)$ 
with $\tau(k(n)) \sim n$ and $k(n) \sim n/A$, it follows that $\sum_{\ell\in\N} \ell\nu(\ell)
=A$. Abbreviate
\be{vnusnu}
v_\nu = \frac{1}{A} \sum_{\ell\in\N} \nu(\ell)\,\E_\mu(Z_\ell),
\quad \sigma^2_\nu =  \frac{1}{A} \sum_{\ell\in\N} \nu(\ell)\,\var_\mu(Z_\ell).
\ee
Write $w-\lim_{n\to\infty}$ to denote convergence in distribution. 

\bt{NCoolLLN+CLT} {\bf [No cooling: Strong LLN and CLT]}
In regime {\rm (R1)} the following hold.\\
{\rm (1)} 
Strong law of large numbers:
\begin{equation}
\lim_{n\to\infty} \frac{X_n}{n} = v_\nu \quad \PP_\QQ\text{-a.s.}
\end{equation}
{\rm (2)}
Central limit theorem:
\be{CLTnocooling}
w-\lim_{n\to\infty} \frac{X_n-v_\nu n}{\sigma_\nu\sqrt{n}} 
= \cN(0,1) \quad \text{ under the law } \PP_\QQ,
\ee
provided
\be{extraass}
\begin{array}{lll}
&{\rm (i)}	 &\bar T^n = o(\sqrt{n}),\\[0.2cm]
&{\rm (ii)}  &\sum_{\ell \in \N} \ell\, |L_n(\ell)-\nu(\ell)| = o\left(1/\sqrt{n}\,\right).
\end{array}
\ee
\et

\noindent
(It is possible to weaken~\eqref{extraass}, but we will not pursue this further.)

\medskip\noindent
$\bullet$ {\bf Cooling.}
Regimes (R2) and (R3) is essence correspond to the situation where the increments of the 
resampling times diverge, i.e., $\lim_{k\to\infty} T_k = \infty$. We have a weak LLN under 
the latter condition only, which we refer to as \emph{cooling}. 

\bt{FCoolLLN} {\bf [Cooling: Weak LLN]}\\
Let $\alpha$ be as in Proposition~{\rm \ref{Solomon}}. If the cooling rule $\tau$ is such that
\be{diverge}
\lim_{k\to\infty}T_k = \infty,
\ee
then
\begin{equation}
w-\lim_{n\to\infty} \frac{X_n}{n} = v_{\mu} \quad \text{ under the law } \PP_\QQ.
\end{equation}
\et

For regimes (R2) and (R3) we derive Gaussian fluctuations for recurrent RWRE:
 
\bt{RecCoolScal} {\bf [Slow and fast cooling: Gaussian fluctuations for recurrent RWRE]}\\
Let $\alpha$ be as in Proposition~{\rm \ref{RecScaling}}. In regimes {\rm (R2)}
and {\rm (R3)},
\begin{equation}
\label{GaussianRec}
w-\lim_{n\to\infty} \frac{X_n-\mathbb{E}_{\mathbb{Q}}(X_n)}{\sqrt{\chi_n(\tau)}} 
= \cN(0,1) \quad \text{ under the law } \PP_\QQ
\end{equation}
with 
\be{RecScales}
\chi_n(\tau)= 
\left\{\begin{array}{ll}
(\sigma^2_\mu \sigma_V)^2 \big(\frac{\beta-1}{\beta}\big)^4 \big(\frac{n}{B}\big)^{1/\beta}\log^4 n, 
&\mbox{ in regime {\rm (R2)}},\\[0.2cm]
(\sigma^2_\mu \sigma_V)^2 \big(\frac{1}{5C^5}\big)\log^5 n, &\mbox{ in regime {\rm (R3)}},
\end{array}
\right. 
\ee
with $\sigma^2_\mu$ the variance of the random variable $\log\rho(0)$ in \eqref{rhodef} 
and $\sigma^2_V \in (0,\infty)$ is the variance of the random variable $V$ in \eqref{densityofV}.
\et

\noindent
Note that the scaling in~\eqref{RecScales} depends on the parameters $(B,\beta)$ 
and $C$ in the regimes (R2) and (R3), respectively, as well as on the law $\mu$ of 
the static random environment.

%%%%%

\subsection{Discussion, open problems and a conjecture}
\label{discussion}

{\bf Preliminary results.}
The results presented in Section~\ref{results} are modest and are only a first step 
in the direction of understanding the effect of the cooling of the random environment.
The weak LLN in Theorem~\ref{FCoolLLN} holds as soon as the cooling is \emph{effective}, 
i.e., the increments of the resampling times diverge (as assumed in \eqref{diverge}). In 
particular, the asymptotic speed $v_\mu$ is the same as for RWRE. This is markedly 
different from Theorem~\ref{NCoolLLN+CLT}, where homogenisation occurs (defined in 
Section~\ref{ssintro}), but with an averaged speed $v_\nu$ that is different from $v_\mu$. 
The CLT in Theorem~\ref{RecCoolScal} can be extended to more general cooling regimes 
than (R2) and (R3), but fails when the cooling becomes too rapid. 

\medskip\noindent
{\bf Future targets.}
In future work we will show that also the strong LLN holds in the effective cooling regime. 
The derivation is more technical and requires that we distinguish between different choices 
of the parameter $s$ defined in \eqref{sparamdef}, which controls the scaling behaviour of 
RWRE (recall Proposition~\ref{TransScaling}). Below we discuss what scaling to expect for 
RWCRE and what properties of the static model in Section~\ref{ssRWRE} would be needed 
to prove this scaling. Essentially, we need \emph{rate of convergence} properties of RWRE, 
as well as control on the \emph{fluctuations of the resampling times}.   
 
\medskip\noindent
{\bf Scaling limits in the recurrent case.}
Suppose that $\alpha$ is as in Proposition~\ref{RecScaling}. Theorem~\ref{RecCoolScal} 
establishes Gaussian fluctuations for the position of the random walk around its mean, of 
an order that depends on the cooling rule $\tau$, given by~\eqref{RecScales}. From our 
knowledge of the static model (recall Proposition~\ref{RecScaling}), we can only conclude 
that $\E_\mu(Z_n) = o(\log^2 n)$. Suppose that, under suitable conditions on $\alpha$, 
we could show that $\E_\mu(Z_n)=o(n^{-1/2} \log n)$ (one example is when $\alpha$ is 
symmetric w.r.t.\ $\tfrac12$, in which case $\E_\mu(Z_n) = 0$ for all $n\in\N_0$). Then, as 
we will see in Section~\ref{s4}, in the slow cooling regime (R2) this extra information would imply for 
RWCRE that $\E_\QQ(X_n)=o(n^{1/2\beta}\log^2 n)$, in which case Theorem~\ref{RecCoolScal} 
would say that $X_n/n^{1/2\beta}\log^2 n$ converges in distribution to a Gaussian random 
variable. In other words, the cooling rule would have the effect of strongly \emph{homogenising} 
the random environment, and the limiting Kesten distribution in~\eqref{densityofV} would 
be washed out. For the recurrent case it is reasonable to expect the presence of a 
\emph{crossover} from a Gaussian distribution to the Kesten distribution as the cooling 
gets slower.

A similar picture should hold in the fast cooling regime (R3). In fact, it would be natural to 
look at even faster cooling regimes, namely, \emph{super-exponential cooling}, in order to 
see whether, after an appropriate scaling, the limiting distribution is the same as in the 
static model. For double-exponential cooling like $\tau(k)=e^{e^k}$ the Lyapunov condition 
in Lemma~\ref{Lindeberg}, on which the proof of Theorem~\ref{RecCoolScal} is based, is 
no longer satisfied.

\medskip\noindent
{\bf Scaling limits in the transient case.}
It is natural to expect that a similar \emph{crossover} also appears in the transient case, 
from the Gaussian distribution to stable law distributions. The general philosophy is the 
same as in the recurrent case: the faster the cooling, the closer RWCRE is to the static 
model and therefore the weaker the homogenisation.

We conjecture the following scenario for the scaling limits of the \emph{centred} 
position of the random walk (i.e., $X_n-v_\mu n$, with $v_\mu$ the speed in 
Theorem~\ref{FCoolLLN}):

%%%%%%%%%%%%%%%%%%%%%%%%%%%%%%%%%%%%%%%%%%%
\begin{table}[htbp]
\vspace{0.5cm}
\centering
\begin{tabular}{c||c|c|c}
\text{Transient Scaling} &(R1) &(R2)  &(R3)\\
\hline
&& $\exists\,\beta_c=\beta_c(s)$: & \\
$s \in (0,2)$ &\emph{CLT} &$\beta<\beta_c$: \emph{Homogenisation} & \emph{Static Law} \\  
&&\hspace{-1.1cm} $\beta>\beta_c$: \emph{Static Law} &\\    
\hline
$s \in (2,\infty)$ &\emph{CLT} &\emph{CLT} &\emph{CLT}\\
\end{tabular}
\end{table}
\vspace{0.5cm}
%%%%%%%%%%%%%%%%%%%%%%%%%%%%%%%%%%%%%%%%%%%%%%%

\noindent
In this table, $s \in (0,\infty)$ is the parameter in Proposition~\ref{TransScaling}, $\beta$ is 
the exponent in regime (R2), and:
\begin{itemize}
\item
\emph{CLT} means that the centred position divided by $\sqrt{n}$ converges in 
distribution to a Gaussian.
\item 
\emph{Homogenisation} means that the centred position divided by a factor that is different 
from $\sqrt{n}$ (and depends on the cooling rule $\tau$) convergences in distribution to a 
Gaussian (compare with Theorem~\ref{RecCoolScal}). 
\item
\emph{Static Law} means that the centred position divided by a factor that is different from 
$\sqrt{n}$ (and depends on the cooling rule $\tau$) has the same limit distribution as in the 
static model (compare with Proposition~\ref{TransScaling}). 
\end{itemize} 
The items under (R2) and (R3) are conjectured, the items under (R1) are proven (compare 
with Theorem~\ref{NCoolLLN+CLT}).

The most interesting feature in the above table is that, in regime (R2) and for $s \in (0,2)$, 
we conjecture the existence of a critical exponent $\beta_c=\beta_c(s)$ above which 
there is Gaussian behaviour and below which there is stable law behaviour. This is motivated 
by the fact that, for $s \in (0,2)$, fluctuations are of polynomial order in the static model 
(see Proposition~\ref{TransScaling}). Hence when the cooling rule is also of polynomial 
order, as in regime (R2), there is a competition between ``localisation between resamplings" 
and ``homogenisation through resamplings''. 
 
\medskip\noindent
{\bf Role of the static model.} 
To establish the CLTs in the second row of the table, for $s \in (2,\infty)$, a deeper 
understanding of the static model is required. In fact, to prove Gaussian behaviour 
in regimes (R2) and (R3) we may try to check a Lyapunov condition, as we do in
the proof of Theorem~\ref{RecCoolScal} (see Lemma~\ref{Lindeberg} below). However, 
as in the proof of Theorem~\ref{RecCoolScal}, $L^p$ convergence for some $p>2$ 
of $Z_n/\sqrt{n}$ would be needed. As far as we know, this is not available in 
the literature. An alternative approach, which would be natural for all regimes and would make 
rigorous the picture sketched in the above table, is to check convergence of the characteristic 
functions of the random variables $X_n-v_\mu n$ properly scaled. The advantage of this 
approach is that, due to the independence of the summands in~\eqref{XRepresentation}, 
this characteristic function factorises into a product of characteristic functions of the 
centred static model, on different time scales. However, we would need suitable 
rate-of-convergence results for these static characteristic functions, which are also
not available.

In the table we did not include the case $s=2$, nor did we comment on the distinction 
between $s \in (0,1)$, $s=1$ and $s\in(1,2)$, for which the centring appears to be 
delicate (compare with Proposition~\ref{TransScaling}).

\medskip\noindent
{\bf Hitting times.}
To derive the static scaling limits in Proposition~\ref{TransScaling}, Kesten, Kozlov and 
Spitzer~\cite{KKS75} (see Peterson~\cite{Peterson} for later work) first derive the scaling 
limits for the hitting times $\sigma_x = \inf\{n\in\N_0\colon\,Z_n=x\}$, $x \in \N$, and afterwards 
pass to the scaling limits for the positions $Z_n$, $n \in \N_0$ via a simple inversion argument. 
This suggests that a further approach might be to look at the hitting times associated with our 
RWCRE model. However, a decomposition into a sum of independent random variables, in 
the flavour of~\eqref{XRepresentation}, would no longer hold for $\sigma_x$, $x \in \N_0$. 
Consequently, in this approach it seems even more difficult to exploit what is known about 
the static model.

%%%%% SECTION 2 %%%%%%%%%%%%%%%%%%%

\section{No cooling: proof of Theorem~\ref{NCoolLLN+CLT}}
\label{s2}

Throughout the sequel we use the same symbol $\PP$ for the annealed law $\PP_\mu$ 
of RWRE in~\eqref{RWREann} and the annealed law $\PP_\QQ$ of RWCRE in 
\eqref{RWCREann}. At all times it will be clear which of the two applies: no confusion 
is possible because they are linked via~\eqref{XRepresentation}.  

\begin{proof}
The proof uses the Lyapunov condition in Lemma~\ref{Lindeberg}.\\ 
(1) Rewrite~\eqref{XRepresentation} as
\be{Beethoven1}
\begin{aligned}
\frac{X_n}{n} &= \frac{1}{n} \sum_{\ell\in\N} \sum_{k=1}^{k(n)} 
Y_k\,{\mathbbm 1}_{\{\tau(k)-\tau(k-1)=\ell\}} + \frac{1}{n}\,\bar  Y^n\\
&\,\hat{=}\, \frac{k(n)}{n}\,\sum_{\ell\in\N} L_n(\ell) 
\left(\frac{1}{k(n) L_n(\ell)}
\sum_{m=1}^{k(n) L_n(\ell)}  Z_\ell^{(m)} \right) + \frac{1}{n}\,Z_{\bar T^n},
\end{aligned}
\end{equation}
where $Z_\ell^{(m)}$, $m\in\N$, are independent copies of $Z_\ell$ and $\hat{=}$ denotes
equality in distribution. Since $\tau(k) \sim Ak$, we have $k(n) \sim n/A$ and $\bar T^n 
= o(n)$. Moreover, since $\abs{Z_m} \leq m$ for all $m \in \N_0$, we have
\begin{equation}
\lim_{n\to\infty} \frac{1}{n}\,Z_{\bar T^n} = 0 \quad \PP\text{-a.s.}, \qquad 
\lim_{N\to\infty} \frac{1}{N} \sum_{m=1}^N Z_\ell^{(m)} = \E(Z_\ell) \quad \PP\text{-a.s.},
\end{equation}
and hence the claim follows by dominated convergence.

\medskip\noindent
(2) Lemma~\ref{Lindeberg} implies that, subject to the Lyapunov condition,
\be{CLTnocooling2}
w-\lim_{n\to\infty} \frac{X_n-a_n}{\sqrt{b_n}} = \mathcal{N}(0,1)
\quad \text{ under the law } \PP
\ee 
with 
\be{anbnformula2}
a_n = \sum_{k=1}^{k(n)} \E(Y_k) + \E(\bar  Y^n), \qquad 
b_n = \sum_{k=1}^{k(n)} \var(Y_k) + \var(\bar  Y^n).
\ee
We need to show that $a_n = v_\nu n + o(\sqrt{n})$ and $b_n \sim \sigma^2_\nu n$,
and verify the Lyapunov condition. 

To compute $b_n$, we note that $\var(\bar  Y^n) \leq (\bar T^n)^2 = o(n)$ by Assumption (i)
in~\eqref{extraass}, and we write 
\be{Bruch1}
\begin{aligned}
\sum_{k=1}^{k(n)} \var(Y_k) 
&= k(n) \sum_{\ell \in \N} L_n(\ell)\,\var(Z_\ell) \\
&= k(n) \sum_{\ell \in \N} \nu(\ell)\,\var(Z_\ell) 
+ k(n) \sum_{\ell \in \N} [L_n(\ell)-\nu(\ell)]\,\var(Z_\ell).
\end{aligned}
\ee
The first term in the right-hand side equals $k(n)\,A\sigma^2_\nu \sim \sigma^2_\nu n$, 
which is the square of the denominator in (\ref{CLTnocooling}). The second term is $k(n)\, 
o(1) = o(n)$, because $\var(Z_\ell) \leq \ell$ and $L_n$ converges to $\nu$ in mean. Hence 
$b_n \sim \sigma^2_\nu n$.
 
To compute $a_n$, we note that $|\E(\bar  Y^n)| \leq \bar T^n = o(\sqrt{n})$ by Assumption 
(i) in \eqref{extraass}, and we write 
\be{Bruch2}
\begin{aligned}
\sum_{k=1}^{k(n)} \E(Y_k) 
&= k(n) \sum_{\ell \in \N} L_n(\ell)\,\E(Z_\ell) \\
&= k(n) \sum_{\ell \in \N} \nu(\ell)\,\E(Z_\ell)
+ k(n) \sum_{\ell \in \N} [L_n(\ell)-\nu(\ell)]\,\E(Z_\ell).
\end{aligned}
\ee
The first term in the right-hand side equals $k(n)\,Av_\nu \sim v_\nu n$, which is the 
numerator in (\ref{CLTnocooling}). By Assumption (ii) in~\eqref{extraass}, the second term is 
$k(n)\,o(1/\sqrt{n})$ because $|\E(Z_\ell)| \leq \ell$. Hence $a_n = v_\nu n + o(\sqrt{n})$.

It remains to verify the Lyapunov condition. To do so, we first show that
\be{maximum}
\max_{1 \leq k \leq k(n)} T_k = o(\sqrt{n}).
\ee 
Define 
\be{defkstar}
M_n = \argmax_{1 \leq k \leq k(n)} T_k,
\ee
i.e., the index for which the gap between two successive resampling times is maximal. If 
there are several such indices, then we pick the largest one. We estimate
\be{gapest}
\begin{aligned}
&\frac{\max_{1 \leq k \leq k(n)}T_k}{\sqrt{n}}
\leq \frac{\tau(M_n) - \tau(M_n-1)}{\sqrt{\tau(M_n)}}\\
&\qquad \leq \frac{1+\big[\big(\tau(M_n)-1\big)
-\tau\big(k\big(\tau(M_n)-1\big)\big)\big]}{\sqrt{\tau(M_n)-1}}
= \frac{1+\bar T^{\tau(M_n)-1}}{\sqrt{\tau(M_n)-1}},
\end{aligned}
\ee
where in the first inequality we use that $\tau(M_n) \leq n$ and in the second inequality 
that $\tau(M_n-1) = \tau(k(\tau(M_n))-1)=\tau(k(\tau(M_n)-1))$ because $m=k(\tau(m))$, 
$m\in\N_0$ (recall Fig.~\ref{fig:resampkn}). By Assumption (i) in~\eqref{extraass}, the 
right-hand side of \eqref{gapest} tends to zero when $\lim_{n\to\infty} M_n = \infty$. If the 
latter fails, then $n \mapsto \max_{1 \leq k \leq k(n)} T_k$ is bounded, and the left-hand 
side of \eqref{gapest} still tends to zero. 

Armed with \eqref{maximum} we prove the Lyapunov condition as follows. Abbreviate 
$\ell(n) = \bar T^n \vee \max_{1 \leq k \leq k(n)} T_k = o(\sqrt{n})$. We have $|Y_k| \leq 
\ell(n)$, $1 \leq k \leq k(n)$, and $|\bar  Y^n| \leq \ell(n)$. This enables us to estimate
\be{sumgapest}
\begin{aligned}
&\frac{1}{b_n^{(2+\delta)/2}} \left[ \sum_{k=1}^{k(n)} 
\E\left(|Y_k - \E(Y_k)|^{2+\delta}\right)  
+ \E\left(|\bar  Y^n - \E(\bar  Y^n)|^{2+\delta}\right)\right]\\
&\leq \frac{[2\ell(n)]^\delta}{b_n^{(2+\delta)/2}}
\left[ \sum_{k=1}^{k(n)} \E\left(|Y_k - E(Y_k)|^2\right)  
+ \E\left(|\bar  Y^n - \E(\bar  Y^n)|^2\right)\right]
= \left(\frac{2\ell(n)}{\sqrt{b_n}}\right)^\delta,   
\end{aligned}
\ee
which tends to zero because $b_n \sim \sigma_\nu^2 n$ and $\ell(n) = o(\sqrt{n})$.
\end{proof}

%%%%% SECTION 3 %%%%%%%%%%%

\section{Cooling: proof of Theorem~\ref{FCoolLLN}}
\label{s3}

\begin{proof}
We need to show that 
\begin{equation}
\lim_{n\to\infty} \PP\left( \left| n^{-1}X_n-v_\mu \right| >\epsilon\right) = 0
\qquad \forall\,\epsilon>0.
\end{equation} 
To that end, we rewrite the decomposition in \eqref{XRepresentation} 
(recall \eqref{taudefs}--\eqref{barTdef}) as
\begin{equation}
\label{Xdec2}
X_n = \sum_{k \in \N} \gamma_{k,n}\,\frac{Y_k}{T_k} + \gamma_n\,\frac{\bar Y^n}{\bar T^n}
\end{equation}
with (see Fig.~\ref{fig:resamp2})
\begin{equation}
\label{alphaGamma}
\gamma_{k,n} = \frac{T_k}{n}{\mathbbm 1}_{\{1 \leq k \leq k(n)\}}, \qquad 
\gamma_n = \frac{\bar T^n}{n}. 
\end{equation}
Note that  $\sum_{k \in \N} \gamma_{k,n} + \gamma_n = 1$ and $\lim_{n\to\infty} 
\gamma_{k,n} = 0$ for all $k \in \N$. 

%%%%%%%%%%%%%%%%%%%%%%%%%%%%%%%%%%%%%%%%%
\begin{figure}[htbp]
\begin{center}
\setlength{\unitlength}{1cm}
\begin{picture}(13,2)(-.5,-1)
\put(.1,.1){\line(1,0){3}}
\put(4.1,.1){\line(1,0){4}}
\put(9.1,.1){\line(1,0){4}}
{\thicklines
\qbezier(13.1,-.1)(13.1,0)(13.1,.25)
}
\qbezier[10](3.4,.1)(3.6,.1)(3.8,.1)
\qbezier[10](8.4,.1)(8.6,.1)(8.8,.1)
\put(0,0){$\bullet$}
\put(2,0){$\bullet$}
\put(5,0){$\bullet$}
\put(7,0){$\bullet$}
\put(10,0){$\bullet$}
\put(-.3,-.5){$\tau(0)$}
\put(1.7,-.5){$\tau(1)$}
\put(4.3,-.5){$\tau(k-1)$}
\put(6.7,-.5){$\tau(k)$}
\put(9.4,-.5){$\tau(k(n))$}
\put(13,-.5){$n$}
\put(5.9,.5){$T_k$}
\put(11.5,.5){$\bar T^n$}
\end{picture}
\end{center}
\vspace{-.5cm}
\caption{Resampling times and increments (recall Fig.~\ref{fig:resamp1}).}
\label{fig:resamp2}
\end{figure}
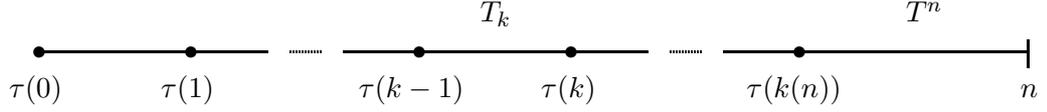
%%%%%%%%%%%%%%%%%%%%%%%%%%%%%%%%%%%%%%%%%

To deal with the representation of $X_n$ in \eqref{Xdec2}, we need the following variant of
a Toeplitz lemma, adapted  to our problem, the proof of which is given in Appendix~\ref{app:Toep}. 

\bl{ToepVariant}
Let $(\gamma_{k,n})_{k,n\in\N}$ and $(\gamma_n)_{n\in\N}$ be as in \eqref{alphaGamma}.
Let $(z_k)_{k\in\N}$ be a real-valued sequence such that $\lim_{k\to\infty} z_k = z^*$ for 
some $ z^*\in\R$. Then
\begin{equation}
\lim_{n\to \infty} \left(\,\sum_{k \in \N} \gamma_{k,n} z_k + \gamma_n z_{\bar T^n}\right) = z^*.
\end{equation}
\el

With the help of \eqref{Xdec2}, we may write
\begin{equation}
\label{go}
n^{-1}X_n-v_\mu =  \left[\,\sum_{k \in \N} \gamma_{k,n} C_k + \gamma_n \bar C^n\right] 
+ \left[\,\sum_{k \in \N} \gamma_{k,n} \left(v_k - v_\mu\right) 
+ \gamma_n \left(\bar v^n - v_\mu\right)\right],
\end{equation}
with 
\begin{equation}
C_k = \frac{Y_k - \E[Y_k]}{T_k},\qquad \bar C^n = \frac{\bar Y^n - \E[\bar Y^n]}{\bar T^n}, 
\qquad v_k = \frac{\E[Z_{T_k}]}{T_k}, \qquad \bar v^n = \frac{\E[Z_{\bar T^n}]}{\bar T^n}.
\end{equation}
Next, note that $\lim_{k\to\infty} v_k = v_\mu$ because $\lim_{k\to\infty} T_k = \infty$. Applying 
Lemma~\ref{ToepVariant} with $z_k=v_k-v_\mu$ and $z^*=0$, we see that the second term 
between square brackets in \eqref{go} tends to zero. Therefore, it suffices to show that the 
first term between square brackets in \eqref{go} tends to zero in probability. 

Estimate
\begin{equation}
\label{bound}
\PP\left(\left|\sum_{k \in \N} \gamma_{k,n} C_k +  \gamma_n \bar C^n \right| > \epsilon\right)
\leq \frac{1}{\epsilon} \left( \sum_{k \in \N} \gamma_{k,n}  \E[|C_k|] +  \gamma_n \E[|\bar C^n|]\right). 
\end{equation}
Since $|C_k| \leq 2$, $k\in\N$, we have $\E[|C_k|]  \leq 2\,\PP(|C_k| >\delta)+\delta$, $k\in\N$, 
for any $\delta>0$. On the other hand, with the help of Proposition~\ref{Solomon} we get
\begin{equation} 
\lim_{k\to\infty} \PP(|C_k| >\delta) = \lim_{k\to\infty} 
\PP\left(\left|\frac{Z_{T_k}-\E[Z_{T_k}]}{T_k}\right| >\delta\right) = 0, 
\qquad \delta>0,
\end{equation}
and hence 
\begin{equation}
\label{zero}
\lim_{k\to\infty}  \E[ |C_k|] = 0.
\end{equation}
Applying Lemma~\ref{ToepVariant} with $z_k= \E[|C_k|]$ and $z^*=0$, and using \eqref{zero}, 
we see that the right-hand side of \eqref{bound} vanishes as $n \to \infty$ for any $\epsilon>0$. 
\end{proof}

%%%%%%% SECTION 4 %%%%%%%%%

\section{Slow and fast cooling: proof of Theorem~\ref{RecCoolScal}}
\label{s4}

The proof again uses the Lyapunov condition in Lemma~\ref{Lindeberg}.

\begin{proof}  
For an arbitrary cooling rule $\tau$, set
\begin{equation}
\chi_n(\tau) = \sum_{k=1}^{k(n)} \var(Y_k) + \var(\bar  Y^n),
\end{equation}
and 
\begin{equation}
\chi_n(\tau;p) = \sum_{k=1}^{k(n)} \E(|Y_k -\E(Y_k)|^p ) 
+ \E(|\bar  Y^n -\E(\bar  Y^n)|^p ), \quad p>2.
\end{equation}
In view of Lemma~\ref{Lindeberg}, it suffices to show that, in regimes (R2) and (R3),
\begin{equation}
\label{criterion}
\lim_{n\to\infty} \frac{\chi_n(\tau;p)}{\chi_n(\tau)^{p/2}} = 0, \qquad p>2.  
\end{equation}

By Proposition~\ref{RecScaling}, we have $\E(Z_n)= o(\log^2 n)$, $n\to\infty$. By 
Proposition~\ref{Lp} below we further have
\begin{equation}
\label{asym}
\E(|Z_n -\E(Z_n)|^2) \sim \Sigma^2\log^4 n, \quad \E(|Z_n -\E(Z_n)|^p )=O(\log ^{2p} n), 
\quad p > 2,
\end{equation}
where $\Sigma=\sigma_\mu^2\sigma_V$. Consequently, using~\eqref{XRepresentation} 
and~\eqref{asym} we get that, for an arbitrary cooling rule $\tau$,
\begin{equation}
\begin{aligned}
\label{scale}
\chi_n(\tau) &\sim \Sigma^2 \sum_{k=1}^{k(n)} \log^4 T_k + \Sigma^2\log^4 \bar T^n,\\
\chi_n(\tau;p) &= \sum_{k=1}^{k(n)} O(\log^{2p} T_k) + O\big(\log^{2p} \bar T^n\big), \quad p>2.
\end{aligned}
\end{equation}
By~\eqref{Tnreg}, in regime (R2) we have $T_k \sim \beta Bk^{\beta-1}$, $k \to \infty$, and 
$\bar T^n \sim (n/B)^{1/\beta}$, $n\to\infty$, while in regime (R3), we have $\log T_k 
\sim Ck$, $k\to\infty$, and $\bar T^n \sim (1/C) \log n$, $n\to\infty$. Thus we see that, 
in both regimes,~\eqref{criterion} holds and $\chi_n(\tau)$ scales as in~\eqref{RecScales}. 
Hence the claim follows from Lemma~\ref{Lindeberg} below.
\end{proof}

%%%%%%% APPENDICES %%%%%%%%%%%

\appendix

%%%%%%%%%%% APPENDIX TOEPLITZ %%%%%%%%

\section{Toeplitz lemma}
\label{app:Toep}

In this appendix we prove Lemma~\ref{ToepVariant}.

\begin{proof}
Estimate
\begin{align}
\label{vai}
\abs{\sum_{k \in \N} \gamma_{k,n} z_k + \gamma_n z_{\bar T^n} - z^*}
\leq \sum_{k \in \N} \gamma_{k,n} \abs{z_k-z^*} + \gamma_n \abs{z_{\bar T^n} - z^*}.
\end{align}
Hence it suffices to shows that, for $n$ large enough, the right-hand side is smaller than 
an arbitrary $\epsilon>0$. To this end, pick $\epsilon_1>0$ (which will be fixed at the end), 
and choose $N_0=N_0(\epsilon_1)$ such that $\abs{z_k-z^*} < \epsilon_1$ for $k>N_0$.
Note further that, in view of \eqref{alphaGamma}, we can choose $N_1=N_1(\epsilon_1)$ 
such that, for $n > N_1$,
\begin{equation}
\label{one}
\sum_{k=1}^{N_0} \gamma_{k,n} \abs{z_k-z^*} < \epsilon_1, 
\end{equation}
and 
\begin{equation}
\label{two}
\gamma_n = \frac{\bar T^n}{n} > \epsilon_1\quad \Longrightarrow \quad  \bar T^n > N_0.
\end{equation}
Write the right-hand side of \eqref{vai} as 
\begin{equation}
\label{sumthree}
\sum_{k=1}^{N_0} \gamma_{k,n} \abs{z_k-z^*} + \sum_{k>N_0} \gamma_{k,n} \abs{z_k-z^*}
+\gamma_n \abs{z_{\bar T^n} - z^*}. 
\end{equation}
The first two terms in \eqref{sumthree} are bounded from above by $\epsilon_1$ for $n>\max\{N_0,N_1\}$. 
Indeed, for the first term this is due to \eqref{one}, while for the second term it is true because 
$\sum_{k \in \N} \gamma_{k,n} \leq 1$ and $\abs{z_k-z^*}<\epsilon_1$ for $k>N_0$. For the 
third term we note that either $\gamma_n=\frac{\bar T^n}{n}>\epsilon_1$, in which case \eqref{two} 
together with $\gamma_n \leq 1$ guarantees that $\gamma_n \abs{z_{\bar T^n} - z^*}<\epsilon_1$, 
or $\gamma_n\leq\epsilon_1$, in which case $\gamma_n\abs{z_{\bar T^n} - z^*}<K\epsilon_1$ with 
$K = \sup_{k \in \N} \abs{z_k-z^*}< \infty$. We conclude that \eqref{sumthree} is bounded from 
above by $(2+\max\{1,K\})\epsilon_1$. Now we let $\epsilon_1$ be such that $(2+\max\{1,K\})\epsilon_1
<\epsilon$, to get the claim. 
\end{proof}

%%%%%% APPENDIX A %%%%%%%%%%%%

\section{Central limit theorem}
\label{app:LL}

\bl{Lindeberg} {\rm ({\bf Lindeberg and Lyapunov condition,} Petrov~\cite[Theorem 22]{P75})}.
Let $U=(U_k)_{k \in \N}$ be a sequence of independent random variables (at least 
one of which has a non-degenerate distribution). Let $m_k=\E(U_k)$ and $\sigma_k^2
=\var(U_k)$. Define
\begin{equation}
\chi_n = \sum_{k=1}^n \sigma_k^2. 
\end{equation}
Then the Lindeberg condition
\begin{equation}
\lim_{n\to\infty} \frac{1}{\chi_n} \sum_{k=1}^n \E\left((U_k -m_k)^2
\,1_{\abs{U_k-m_k} \geq \epsilon \sqrt{\chi_n}}\right) = 0
\end{equation}
implies that
\begin{equation}
w-\lim_{n\to\infty} \frac{1}{\sqrt{\chi_n}} \sum_{k=1}^n (U_k-m_k) = \cN(0,1).
\end{equation}
Moreover, the Lindeberg condition is implied by the Lyapunov condition
\be{Lyapunov}
\lim_{n\to\infty} \frac{1}{\chi_n^{p/2}} \sum_{k=1}^n 
\E\left(|U_k -m_k|^p\right) = 0 \quad \text{ for some } p > 2.
\ee
\el

%%%% APPENDIX B %%%%%%%%%%%%

\section{RWRE: $L^p$ convergence under recurrence}
\label{appLp}

The authors are grateful to Zhan Shi for suggesting the proof of Proposition~\ref{Lp}
below.

We begin by observing that all the moments of the limiting random variable $V$ in 
Theorem~\ref{RecScaling} are finite.

\bl{pMom} 
Let $V$ be the random variable with density function~\eqref{densityofV}. Let $P$ denote
its law. Then $E\left(V^p\right)<\infty$ for all $p>0$ with $E\left(V^{2k}\right)=0$ for 
$k\in\mathbb{N}$.
\el

\begin{proof}
For $k\in\mathbb{N}$, it follows from~\eqref{densityofV} that $E(V^{2k})=0$. 
For arbitrary $p>0$, compute
\be{pthmoment}
E\left(|V|^p\right)  
= \frac{4}{\pi}\sum_{k\in\N_0} \frac{(-1)^k}{2k+1}
\int_0^\infty x^p\,\exp\left(-\frac{(2k+1)^2\pi^2}{8}x \right) \dd x.
\ee
Since $b^q \int_0^\infty x^{q-1}\,e^{-bx}\,\dd x = \Gamma(q)$, the integral in 
\eqref{pthmoment} equals
\be{gamma3}
\frac{8^{p+1}\Gamma(p+1)}{(2k+1)^{2(p+1)}\pi^{2(p+1)}}.
\ee
Therefore
\be{pthmoment2}
E\left(|V|^p\right) = \frac{4\Gamma(p+1)8^{p+1}}{\pi^{2p+3}}
\sum_{k\in\N_0} \frac{(-1)^k}{(2k+1)^{2p+3}},
\ee
which is finite for all $p>0$.
\end{proof}

\bp{Lp}
The convergence in Proposition~{RecScaling} holds in $L^p$ for all $p>0$. 
\ep

\begin{proof}
The proof comes in 3 Steps.

\medskip\noindent
{\bf 1.}
As shown by Sinai~\cite{S82},
\begin{equation}
\label{cvg_Sinai}
w-\lim_{n\to\infty} \frac{Z_n - b_n}{\log^2 n} = 0 \quad \text{under the law } \PP,
\end{equation}
where $b_n$ is the bottom of the valley of height $\log n$ containing the origin for the 
potential process $(U(x))_{x\in\Z}$ given by
\begin{equation}
U(x) = \left\{\begin{array}{ll}
\sum_{y=1}^x \log \rho(y), &x \in \N,\\
0, &x=0,\\
- \sum_{y=x}^{-1} \log \rho(y), &x \in -\N,
\end{array}
\right.
\end{equation}
with $\rho(y) = (1-\bar\omega(y))/\bar\omega(y)$. This process depends on the environment
$\omega$ only, and 
\begin{equation}
w-\lim_{n\to\infty} \frac{b_n}{\log^2 n} = V \quad \text{under the law } \alpha^\Z. 
\end{equation}
We will prove the claim by showing that, for all $p>0$, 
\begin{equation}
\label{tightness}
\sup_{n \geq 3} \,\E_{\alpha^\Z} \Big(\Big|\frac{b_n}{\log^2 n}\Big|^p\Big) < \infty,
\qquad 
\sup_{n \geq 3} \,\E\Big(\Big|\frac{Z_n}{\log^2 n}\Big|^p\Big) < \infty.   
\end{equation}
To simplify the proof we may assume that there is a \emph{reflecting barrier at the origin}, 
in which case $b_n$ and $Z_n$ take values in $\N_0$. This restriction is harmless because 
without reflecting barrier we can estimate $|b_n| \leq \max\{b_n^+,-b_n^-\}$ and $|Z_n| 
\leq \max\{Z_n^+,-Z_n^-\}$ in distribution for two independent copies of $b_n$ and $Z_n$ 
with reflecting barrier to the right, respectively, to the left. 

\medskip\noindent
{\bf 2.}
To prove the first half of~\eqref{tightness} with reflecting barrier, define
\begin{equation}
H(r) = \inf\{x \in \N_0\colon\, |U(x)| \geq r\}, \qquad r \geq 0.
\end{equation}
Then
\begin{equation}
\label{ub0}
b_n \leq H(\log n).
\end{equation}
We have
\begin{equation}
\label{ub1}
\E_{\alpha^\Z}\Big(\Big|\frac{H(\log n)}{\log^2 n}\Big|^p\Big) 
= \int_0^\infty p \lambda^{p-1}\,\PP\big(H(\log n) > \lambda \log^2 n\big)\,d\lambda.
\end{equation}
Since $\int_0^1 p \lambda^{p-1} d\lambda = 1$, we need only care about $\lambda 
\geq 1$. To that end, note that
\begin{equation}
\label{ub2}
\big\{H(\log n) > \lambda \log^2 n\big\} 
= \Big\{\max_{0 \leq x \leq \lambda \log^2 n} |U(x)| < \log n\Big\}
\end{equation} 
and
\begin{equation} 
\label{ub3}
\alpha^\Z\left(\max_{0 \leq x \leq \lambda \log^2 n} |U(x)| < \log n\right)
\leq \hat{P}\left(\max_{0 \leq t \leq \lambda N} \sigma |W(t)| < \sqrt{N}\,\right),
\qquad N = \log^2 n, 
\end{equation}
where $\sigma^2$ is the variance of $\rho(0)$ and $(W(t))_{t\geq 0}$ is standard 
Brownian motion on $\R$ with law $\hat{P}$. But there exists a $c>0$ (depending on
$\sigma$) such that
\begin{equation}
\label{ub4}
\hat{P}\Big(\max_{0 \leq t \leq \lambda N} \sigma |W(t)| < \sqrt{N} \Big)
= \hat{P}\Big(\max_{0 \leq t \leq \lambda} \sigma |W(t)| < 1 \Big)
\leq e^{-c\lambda}, \qquad \lambda \geq 1.
\end{equation} 
Combining~\eqref{ub0}--\eqref{ub4}, we get the first half of~\eqref{tightness}.

\medskip\noindent
{\bf 3.}
To prove the second half of~\eqref{tightness}, write
\begin{equation}
\E\Big(\Big|\frac{Z_n}{\log^2 n}\Big|^p\Big)
= \int_0^\infty p \lambda^{p-1}\,
\PP\big(Z_n > \lambda \log^2 n\big) \,d\lambda.
\end{equation}
Again we need only care about $\lambda \geq 1$. As in Step 3, we have
\begin{equation}
\sup_{n\geq 3} \,\int_1^\infty p \lambda^{p-1}\,
\alpha^\Z\big(H(\lambda^{1/3} \log n) > \lambda \log^2 n\big) \,d\lambda < \infty.
\end{equation}
It therefore remains to check that
\begin{equation}
\label{ub2alt}
\sup_{n \geq 3} \, \int_1^\infty p \lambda^{p-1}\,\PP(\cE_{\lambda,n})\,d\lambda < \infty
\end{equation}
with
\begin{equation} 
\cE_{\lambda,n} = \Big\{Z_n > \lambda \log^2 n,\,
H(\lambda^{1/3} \log n) \leq \lambda \log^2 n\Big\}.
\end{equation}
To that end, for $x\in \N_0$, let $T(x) = \inf\{n\in\N_0\colon\,Z_n =x\}$. On the event 
$\cE_{\lambda,n}$ we have $T(H(\lambda^{1/3} \log n)) \leq n$. Therefore, by 
Golosov~\cite[Lemma 7]{G86},
\begin{equation}
\PP(T(x) \leq n \mid \omega) \leq n \,\exp\Big(- \max_{0 \leq y<x} [U(x-1)-U(y)]\Big),
\qquad x \in \N,\,n \in \N,
\end{equation}
which is bounded from above by $n\,e^{-U(x-1)}$. Picking $x = H(\lambda^{1/3} \log n)$, we 
obtain
\begin{equation}
\PP(\cE_{\lambda,n} \mid \omega) \leq n \,e^{-U(H(\lambda^{1/3} \log n)-1)},
\end{equation}
which is approximately $n\,e^{-\lambda^{1/3} \log n}$ because $U(H(x))$ is approximately
$x$. (The undershoot at $x$ can be neglected because it has finite first moment, by our 
assumption that $\sigma<\infty$.) Taking the expectation over $\omega$, we get
\begin{equation}
\PP(\cE_{\lambda,n}) \leq n^{-(\lambda^{1/3}-1)}.
\end{equation}
This implies~\eqref{ub2alt}, and hence we have proved the second half of~\eqref{tightness}.
\end{proof}

%%%%%%% BIBLIOGRAPHY %%%%%%%%%%%%


\begin{thebibliography}{99}

\bibitem{AdGCdCpr}
L.\ Avena, Y.\ Chino, C.\ da Costa and F.\ den Hollander,
Random walk in cooling random environment: strong law of large numbers and 
large deviation principle, in preparation.

\bibitem{BGT87}
N.H.\ Bingham, C.M.\ Goldie and J.L.\ Teugels,
\emph{Regular Variation},
Encyclopedia of Mathematics and its Applications, Vol.\ 27, Cambridge University Press,
Cambridge, 1987. 

\bibitem{BMPseries}
C.\ Boldrighini, R.A.\ Minlos, A.\ Pellegrinotti and E.\ Zhizhina,
Continuous time random walk in dynamic random environment,
Markov Proc.\ Relat.\ Fields 21 (2015) 971--1004.

\bibitem{G86}
A.O.\ Golosov, 
Localization of random walks in one-dimensional random environments, 
Comm.\ Math.\ Phys.\ 92 (1986) 491--506.
  
\bibitem{K86}
H.\ Kesten, 
The limit distribution of Sinai's random walk in random environment,
Physica 138A (1986) 299--309.

\bibitem{KKS75}
H.\ Kesten, M.W.\ Kozlov and F.\ Spitzer, 
A limit law for random walk in random environment,
Compositio Math.\ 30 (1975) 145--168. 

\bibitem{Peterson}
J.\ Peterson, 
\emph{Limiting Distributions and Large Deviations for Random Walks in Random Environments},
PhD thesis, University of Minnesota, July 2008. [arXiv:0810.0257]

\bibitem{P75} 
V.\ Petrov, 
\emph{Sums of Independent Random Variables}, 
Springer, 1975.

\bibitem{RAS05}
F.\ Rassoul-Agha and T.\ Seppalainen, 
An almost sure invariance principle for random walks in a space-time random 
environment, 
Probab.\ Theory Relat.\ Fields 133 (2005) 299--314.

\bibitem{S82}
Ya.G.\ Sinai, 
The limiting behavior of a one-dimensional random walk in a random medium,
Theory Prob.\ Appl.\ 27 (1982) 256--268. 

\bibitem{S75}
F.\ Solomon, 
Random walks in random environment, 
Ann.\ Prob.\ 3 (1975) 1--31.  

\bibitem{Z04}
O.\ Zeitouni, 
Random walks in random environment, 
XXXI Summer School in Probability, Saint-Flour, 2001, 
Lecture Notes in Math.\ 1837 (2004) 189--312.

\end{thebibliography}
\end{document}